\documentclass{amsart}
\usepackage[T2A]{fontenc}
\usepackage[utf8]{inputenc}
\usepackage[english]{babel}
\usepackage{amsmath}
\usepackage{amssymb}
\usepackage{bbm}

\usepackage{xcolor}

\numberwithin{equation}{section}
\newtheorem{theorem}{Theorem}[section]
\newtheorem{lem}[theorem]{Lemma}
\newtheorem{pr}[theorem]{Proposition}

\newtheorem*{theorem*}{Theorem}
\newtheorem*{lem*}{Lemma}
\theoremstyle{definition}

\newtheorem{rem}[theorem]{Remark}
\newtheorem{problem}[theorem]{Problem}

\newcommand{\si}{\sigma}
\newcommand{\sid}{\sigma_d}
\newcommand{\meal}{m}
\newcommand{\za}{\zeta}

\newcommand{\hol}{\mathcal{H}ol}
\newcommand{\Nbb}{\mathbb N}
\newcommand{\Rbb}{\mathbb R}
\newcommand{\Dbb}{\mathbb D}
\newcommand{\DDb}{\overline{\mathbb D}}
\newcommand{\Tbb}{\partial \mathbb{D}}
\newcommand{\T}{\mathbb{T}}

\newcommand{\eps}{\varepsilon}

\newcommand{\bd}{B_d}
\newcommand{\spd}{\partial B_d}
\newcommand{\ccd}{\mathbb{C}^d}
\newcommand{\bbd}{\overline{B}_d}
\newcommand{\MM}{M_+}
\newcommand{\qball}{Q}
\newcommand{\hp}{{H^p}}

\newcommand{\rad}{\mathcal{R}}

\newcommand{\D}{\mathbb{D}}

\begin{document}
	
	\title[Reverse Carleson measures for spaces of analytic functions]{Reverse Carleson measures for spaces of analytic functions}
	
	\author{Evgueni Doubtsov}
	\address{
		St. Petersburg Department of Steklov Mathematical Institute, Fontanka 27, St. Petersburg 191023, Russia}
	\email{dubtsov@pdmi.ras.ru}

	\author{Anton Tselishchev}
	\address{
		St. Petersburg Department of Steklov Mathematical Institute, Fontanka 27, St. Petersburg 191023, Russia}
	\email{celis\_anton@pdmi.ras.ru}
	
	\author{Ioann Vasilyev}
	\address{
		St. Petersburg Department of Steklov Mathematical Institute, Fontanka 27, St. Petersburg 191023, Russia}
	\email{ivasilyev@pdmi.ras.ru}

	\thanks{The work of A.~Tselishchev and I.~Vasilyev was supported by the Ministry of Science and Higher Education of the Russian Federation (agreement 075-15-2025-344 dated 29/04/2025 for Saint Petersburg Leonhard Euler International Mathematical Institute at PDMI RAS)}

	\begin{abstract}
		Let $X$ be a quasi-Banach space of analytic functions in the unit disc and let $q > 0$. 
		A finite positive Borel measure $\mu$ on the unit closed disc $\overline{\mathbb{D}}$
		is called a \emph{$q$-reverse Carleson measure} for $X$ if and only if there exists a constant $C > 0$ such that
		\[
		\|f\|_{X}\leq C \|f\|_{L^q(\overline{\D},d\mu)}
		\]
		for all $f\in X\cap C(\overline{\D})$.
		We fully characterize the $q$-reverse Carleson measures with all $q > 0$ for Hardy spaces $H^p(\Dbb)$ with all $0<p\leq \infty$, for the space $\mathrm{BMOA}(\Dbb)$ and for the Bloch space.
		In addition, we describe $q$-reverse Carleson measures for the holomorphic Triebel--Lizorkin spaces $HF_0^{q,r}$ and the
		holomorphic Besov spaces $HB_0^{q,r}$.
		Related results are obtained for the Hardy spaces and certain holomorphic Triebel--Lizorkin spaces in the unit ball of $\mathbb{C}^d$.
	\end{abstract}
	
	\maketitle
	
	\section{Introduction}\label{s_int}
	
	\subsection{Hardy spaces and reverse Carleson measures}
	
	Let $\mathcal{H}ol(\Dbb)$ denote the space of holomorphic functions in the unit disc $\Dbb$ of $\mathbb{C}$.
	For $0 < p < \infty$, the Hardy space $H^p=H^p(\D)$ consists of $f\in \mathcal{H}ol(\Dbb)$ such that
	$$
	\|f\|_{H^p}^p = \sup_{0 < r < 1} \int_0^1 |f(r e^{2\pi i t})|^p\, dt < \infty.
	$$
	The space $H^\infty = H^\infty(\mathbb D)$ is defined as the vector space of bounded holomorphic functions in the disc, with the norm
	$$
	\|f\|_{H^\infty} = \sup_{|z|<1} |f(z)|.
	$$
	Note that the space $H^p$ is a Banach space for $p \ge 1$ but it is merely a quasi-Banach space if $0 < p < 1$.
	
	Let $\MM(\D)$ denote the space of finite positive Borel measures on $\D$. 
	A classical theorem by L. Carleson \cite{Ca58, MR141789} tells us that the space
	$H^p(\Dbb)$ embeds into $L^p(\Dbb, \mu)$ for $\mu\in \MM(\D)$ meaning that there exists a constant $C > 0$ such that 
	\begin{equation}\label{e_Carl_1}
		\| f \|_{L^p(\Dbb, \mu)} \le C \|f\|_\hp \quad\textrm{for all}\ f \in\hp
	\end{equation}
	if and only if the measure 
	$\mu$ satisfies the following condition of L. Carleson: there exists a constant  $C > 0$ such that 
	\begin{equation}\label{e_Carl_geom_1}
		\mu(S_I) \le C \meal(I)\quad\textrm{for all arcs}\ I\subset\Tbb,
	\end{equation}
	where
	\[
	S_I = \left\{z \in\overline{\Dbb}:\ 1-\meal(I) \le |z| \le 1,\ \frac{z}{|z|} \in  I \right\}
	\]
	is a standard Carleson window. Here $m$ stands for Lebesgue measure on $\partial\D$. We will sometimes write $|I|$ instead of $m(I)$ further in the text.
	
	The topic of Carleson measures, that is, measures satisfying condition \eqref{e_Carl_1}, is nowadays very well studied. In particular, the famous Carleson's theorem has been generalized to various spaces of analytic functions in the unit disc and in other domains.
	
	In the present paper we are interested in the problem of  characterization of \emph{reverse Carleson measures}: that is, roughly speaking, of measures $\mu$ which satisfy the inequality reverse to \eqref{e_Carl_1}. The study of such measures (in the context of Bergman spaces) was started by D. Luecking \cite{Lu85}.
	However, the first characterization of reverse Carleson measures for Hardy spaces $H^p$ with $1 < p < \infty$ was obtained only rather recently  by A. Hartmann, X. Massaneda, A. Nicolau, and J. Ortega-Cerd\`a \cite{HMNO14}. 
	
	Let us now pass to the exact definitions. A measure $\mu\in M_+(\overline{\D})$ is called a reverse Carleson measure for the space $H^p$ if there exists a constant  $C > 0$ such that the  inequality
	\begin{equation}\label{reverse_carl_disc}
		\|f \|_\hp \le C \|f\|_{L^p(\overline{\D},\mu)}
	\end{equation}
	holds for all functions $f\in  C(\overline{\D}) \cap H^p$. Note that here we consider the measures $\mu$ supported in the \emph{closed} unit disc $\overline{\D}$; the reason for that will be clear from what follows.
	It turns out that the \emph{balayage} phenomenon occurs in the problem mentioned in the previous paragraph in a more general context: only the part of the measure which lives on the boundary $\partial \D$ is `responsible' for the reverse Carleson inequality; see Theorem~\ref{t_Carl_Hp_onedim} below.
	
	It is worth noting that the classical theorem about direct Carleson measures for spaces $\hp$ can as well be formulated for measures $\mu\in \MM(\overline{\D})$.
	Indeed, in this case the condition \eqref{e_Carl_1} does make sense for  functions $f$
	that belong to the intersection $C(\overline{\D}) \cap H^p$, which is dense in the space $H^p$.
	The condition \eqref{e_Carl_geom_1} in turn gives that the restriction $\mu |_{\partial\D}$ is
	absolutely continuous with respect to Lebesgue measure $m$ and the Radon--Nikodym derivative ${d\mu |_{\partial\D}}/{d m}$ is bounded.
	
	If a measure $\mu\in\MM(\DDb)$ satisfies condition \eqref{e_Carl_geom_1} (that is, if one supposes a priori that $\mu$ is a direct Carleson measure), P. Lef\`evre, D. Li, H. Queff\'elec, and L. Rodriguez-Piazza~\cite{LLQR12} proved that the property \eqref{reverse_carl_disc} for $1 < p < \infty$ is equivalent to the following condition:
	\[
	\mu(S_I ) \ge C\meal(I)\quad\textrm{for all arcs}\ I \subset \Tbb.
	\]
	Here and everywhere else below  $C$  stands for a positive absolute  constant.
	
	Further, in \cite{HMNO14}, the same was proved for all $\mu\in\MM(\DDb)$. 
	This result was recently generalized to higher dimensions, more specifically, to Hardy spaces in the unit ball of $\mathbb{C}^d$; see \cite{Dou24}.
	
	See also \cite{aguilar2025inequalities} for a recent connected result on reverse Carleson measures in tent spaces.
	
	\subsection{New results and organization of the paper}
	The case of exponents $p$ between $0$ and $1$ was not considered in \cite{HMNO14}, 
	however, the corresponding result can be easily deduced from the description of reverse Carleson measures in the case $1 < p < \infty$, see Theorem~\ref{t_Carl_Hp_onedim} below.
	
	Our paper is (mostly) dedicated to various generalizations of this theorem. Let us briefly describe our main results.
	
	In Section~\ref{s_hp}, we extend the result from~\cite{HMNO14} to all exponents $p > 0$. 
	This generalization is almost immediate, however, we prefer to present here a short self-contained proof, since some of our further results rely on it. Overall, Theorem~\ref{t_Carl_Hp_onedim} below serves as an important starting point for our further considerations.
	
	In Section~\ref{s_pq}, we describe so-called $q$-reverse Carleson measures for the space $H^p$, that is, the measures $\mu\in M_+(\overline{\D})$ such that for every $f\in H^p(\D)\cap C(\overline{\D})$ the inequality 
	$$
	\|f\|_{H^p(\D)}\leq C \|f\|_{L^q(\overline{\D};\mu)}
	$$
	holds with some constant $C > 0$. We obtain a complete answer to this problem for all values $0 < p, q < \infty$. It is important to note that our method here is completely different from the  methods used in~\cite{HMNO14}: in fact, our short proof uses ideas from harmonic (and not complex) analysis. For the results concerning the corresponding problem for direct Carleson measures, see \cite{Dur69, Lu91}.
	
	In Section~\ref{s_Bloch}, we completely classify reverse Carleson measures for the Bloch space. In more detail, we prove that given any $q>0$, there are no finite $q$-reverse Carleson measures for the Bloch space, for the space $\mathrm{BMOA}(\Dbb)$ and for the Hardy space $H^\infty(\Dbb).$ In the proof of the main result of the fourth section we will use certain properties of analytic Fej\'er polynomials. We stress that here as well we completely dispense with the methods used in~\cite{HMNO14}. Notice that the corresponding classification of direct Carleson measures for the Bloch space remains an open problem, see e.g. papers~\cite{girela2008carleson} and~\cite{bao2024carleson}.
	
	Section~\ref{s_Triebel} is dedicated to the description of reverse Carleson measures for the holomorphic Triebel--Lizorkin spaces. The scale of Triebel--Lizorkin spaces includes all Hardy spaces, and we prove that, depending on the parameters, either the analogue of Theorem~\ref{t_Carl_Hp_onedim} holds or there are no reverse Carleson measures for these spaces. In the course of the proof, we will slightly generalize a result of Rudin \cite{Ru55}; see Lemma~\ref{q-var} below.
	
	In Section~\ref{s_besov}, we consider yet another important scale of function spaces: that is, Besov spaces of holomorphic functions. Using the recent results from~\cite{baranov2024analytic}, we prove that for these spaces either the analogue of Theorem~\ref{t_Carl_Hp_onedim} holds or reverse Carleson measures do not exist. Note that the question of describing the reverse Carleson measures for Besov spaces was mentioned at the end of~\cite{hartmann2021dominating}.
	
	In Section~\ref{s_highdim}, we consider the spaces of holomorphic functions in the unit ball of $\mathbb{C}^d$. It turns out that it is much more difficult to obtain a description of reverse Carleson measures in this context, therefore, we consider only multidimensional Hardy and certain Triebel--Lizorkin spaces.
	
	Finally, in Section~\ref{s_remarks} we make some concluding remarks and formulate three open questions about reverse Carleson measures.

	\subsection*{Notations}
	Throughout the rest of the text $C>0$ denotes a harmless constant; the signs $\lesssim$ and $\gtrsim$ indicate that the left-hand (right-hand) side of an inequality is less than the right-hand (left-hand) side multiplied
	by such a harmless constant. 
	The sign $\asymp$ means that we have both bounds  $\lesssim$ and  $\gtrsim$ at the same time.
	
	\subsection*{Acknowledgment}
The authors are kindly grateful to the anonymous referee for many valuable comments and suggestions that helped to improve our manuscript.
	
	\section{Hardy spaces in the unit disc}\label{s_hp}
	
	\subsection{Formulation of the result}
	
	Before we proceed to the description of reverse Carleson measures for Hardy spaces, let us introduce the following notation: for $\lambda\in\D$, we put
	$$
	k_\lambda(z) = \frac{1}{1-\overline{\lambda}z},\quad z\in\D.
	$$
	Let $l$ be an arbitrary integer such that $pl > 1$. It turns out that it is enough to check the inequality \eqref{reverse_carl_disc} for the functions $f = k_\lambda^l$. Denote
	\[
	K_{\lambda,l} = \frac{k_\lambda^l}{\|k_\lambda^l\|_{\hp}}.
	\]
	
	Then the following theorem holds.
	
	\begin{theorem}\label{t_Carl_Hp_onedim}
		Let $0 < p < \infty$ and let $\mu \in \MM(\overline{\D})$. Let $l$ be an arbitrary integer such that $pl > 1$.
		The following assertions are equivalent.
		\begin{enumerate}
			\item[(i)] There exists a constant $C > 0$ such that
			\[
			\int_{\overline{\D}}
			|f|^p\,d\mu \ge C \|f\|_\hp^p\quad\textrm{for all}\ f \in C(\overline{\D}) \cap H^p.
			\]
			\item[(ii)] There exists a constant $C > 0$ such that
			\[
			\int_{\overline{\D}}
			|K_{\lambda,l}|^p\,d\mu \ge C \quad\textrm{for all}\ \lambda\in\D.
			\]
			\item[(iii)] There exists a constant $C > 0$ such that
			\[
			\mu(I) \ge C m(I)\quad\textrm{for all arcs}\ I \subset \partial\D.
			\]
		\end{enumerate}
	\end{theorem}
	
	We note here that the question about the description of reverse Carleson measures for spaces $\hp(\D)$, $0<p<\infty$,
	has no solutions for measures $\mu$ defined in the open unit disc  $\D$.
	Indeed, this fact is a direct consequence of the property (iii) in the above theorem.
	
	In \cite{HMNO14}, Theorem~\ref{t_Carl_Hp_onedim} is proved for $p > 1$ and $l = 1$. Now suppose that $p \le 1$ and $pl > 1$. Clearly, the only nontrivial implication of this theorem is $(\mathrm{ii})\Rightarrow (\mathrm{iii})$. However, condition (ii) can be rewritten as follows:
	$$
	\|k_\lambda^l\|_{L^p(\mu)} \ge c \|k_\lambda^l\|_{H^p},
	$$
	which in turn is equivalent to the estimate
	$$
	\|k_\lambda\|_{L^{pl}(\mu)} \ge c \|k_\lambda\|_{H^{pl}}.
	$$
	Since $pl > 1$, now condition~(iii) follows from \cite[Theorem 2.1]{HMNO14}.
	
	We have just explained how to prove Theorem~\ref{t_Carl_Hp_onedim} in a very short way. However, below we present a short self-contained proof (which is very similar to the argument from~\cite{HMNO14}),
	since we use the technical lemmas from this proof in Section~\ref{s_Triebel} below.
	
	\subsection{Essential computations}
	
	In the proof of Theorem~\ref{t_Carl_Hp_onedim} we will need two simple computational lemmas. The first of them is a standard fact (see e.g. \cite[Proposition 1.4.10]{Ru80}); we omit its proof.
	
	\begin{lem}\label{very_simple}
		For any $q > 1$ and $|\lambda| < 1$, we have
		$$
		\int_0^1 \frac{dt}{|1-\lambda e^{2\pi i t}|^q} \asymp (1-|\lambda|)^{1-q}.
		$$
	\end{lem}
	
	For any arc $I\subset\partial\D$, let $S_{I,h}$ denote the $h$-Carleson window:
	$$
	S_{I,h} = \left\{z\in\overline{\D}: 1-h\le |z|\le 1,\, \frac{z}{|z|}\in I \right\}.
	$$
	
	Let $A$ denote Lebesgue measure on $\mathbb{C}$.

	\begin{lem}\label{simple}
		Let $q > 1$ and let $I\subset\partial\D$ be a closed arc. For $0 < h < 1$, put
		\begin{equation}
			\Phi_h(z) = \frac{1}{h}\int_{S_{I,h}} \frac{(1-|\lambda|)^{q-1}}{|1-\overline{\lambda}z|^q}\, dA(\lambda),\quad z\in\overline{\D}.
		\end{equation}
		Then
		\begin{itemize}
			\item there exists a constant $C=C(q)$ such that $|\Phi_h(z)|\le C$ for all $z\in\overline{\D}$ and $0 < h \le 1$;
			
			\item $\lim_{h\to 0} \Phi_h(z) = 0 $ for all $z\in\overline{\D}\setminus I$.
		\end{itemize}
	\end{lem}
	
	\begin{proof}
		We start with the proof of the first assertion of the lemma. Obviously, we have
		$$
		|\Phi_h(z)|\le \frac{1}{h}\int_{1-h\le|\lambda|\le 1} \frac{(1 - |\lambda|)^{q-1}}{|1-\overline{\lambda}z|^q}\, dA(\lambda).
		$$
		Without loss of generality, we can assume that $z\in\mathbb{R}$ and $z > 0$. Integrating in polar coordinates and using Lemma~\ref{very_simple}, we write
		\begin{multline*}
			|\Phi_h(z)|\le \frac{2\pi}{h} \int_{1-h}^1 \int_0^1 \frac{(1-r)^{q-1} r}{|1-rze^{2\pi i t}|^q }\, dt\, dr \lesssim \frac{1}{h} \int_{1-h}^1 \frac{(1-r)^{q-1} r}{(1-rz)^{q-1}}\, dr \\ \le \frac{1}{h} \int_{1-h}^1 \frac{(1-r)^{q-1} }{(1-rz)^{q-1}}\, dr\le\frac{1}{h} \int_{1-h}^1 \frac{(1-r)^{q-1}}{(1-r)^{q-1}}\, dr = 1.
		\end{multline*}
		
		Now, we prove the second assertion of the lemma. If $z\in\overline{\mathbb{D}}\setminus I$, then there exists a $\delta > 0$ such that $|1-\overline{\lambda}z| \ge\delta$ for every $\lambda\in S_{I,h}$. Hence, we have
		$$
		|\Phi_h(z)| \leq \frac{1}{h\delta^q} \int_{S_{I,h}} h^{q-1}\, dA(\lambda)\leq \frac{1}{\delta^q} |I| h^{q-1} \xrightarrow[h\to 0]{} 0,
		$$
		and Lemma~\ref{simple} follows.
	\end{proof}
	
	\subsection{Proof of Theorem~\ref{t_Carl_Hp_onedim}}
	
	The implication $(\mathrm{i})\Rightarrow (\mathrm{ii})$ is obvious. It is also easy to see that condition (iii) implies (i). Indeed, condition (iii) means that Lebesgue measure $m$ on $\partial\D$ is absolutely continuous with respect to the measure $\mu|_{\partial\D}$ and its density is bounded; clearly, this implies condition (i).
	
	It remains to prove the implication $(\mathrm{ii})\Rightarrow (\mathrm{iii})$. Due to Lemma~\ref{very_simple}, we have 
	\begin{equation}\label{norm_klambda}
		\|k_\lambda^l\|^p_{H^p} \asymp (1-|\lambda|)^{1-pl}.
	\end{equation}
	Therefore, condition (ii) (which in turn follows if we substitute $f=k_\lambda^l$ into condition (i)) is equivalent to the following inequality:
	$$
	\int_{\overline{\D}} \frac{d\mu(z)}{|1-\overline{\lambda}z|^{pl}}\gtrsim (1-|\lambda|)^{1-pl}.
	$$
	We can rewrite this inequality as follows:
	$$
	\int_{\overline{\D}} \frac{(1-|\lambda|)^{pl-1}}{|1-\overline{\lambda}z|^{pl}}\, d\mu(z)\gtrsim 1.
	$$
	Integrating the last inequality over the set $S_{I,h}$, we get
	$$
	\int_{\overline{\D}}\int_{S_{I,h}} \frac{(1-|\lambda|)^{pl-1}}{|1-\overline{\lambda}z|^{pl}}\, dA(\lambda)\, d\mu(z)\gtrsim |I|h.
	$$
	
	Putting 
	$$
	\Phi_h(z) = \frac{1}{h}\int_{S_{I,h}}\frac{(1-|\lambda|)^{pl-1}}{|1-\overline{\lambda}z|^{pl}}\, dA(\lambda),
	$$
	we arrive at the following condition:
	\begin{equation}\label{estimate_on_Ф}
		\int_{\overline{\D}} |\Phi_h(z)|\, d\mu(z)\gtrsim |I|.
	\end{equation}
	Note that the function $\Phi_h$ has the same form as the function from Lemma~\ref{simple}. We can write the integral over $\overline{\D}$ as the sum of integrals over $I$ and over $\overline{\D}\setminus I$. Due to Lemma~\ref{simple}, we can use the dominated convergence theorem in order to conclude that
	$$
	\lim_{h\to 0}\int_{\overline{\D}\setminus I} |\Phi_h(z)|\, d\mu(z) = 0.
	$$
	On the other hand, 
	$$
	\int_I |\Phi_h(z)|\, d\mu(z)\le C\mu(I).
	$$
	Overall, letting $h$ to $0$ in equation \eqref{estimate_on_Ф}, we conclude that $m(I)\le C\mu(I)$, and we are done.

	\section{$(p,q)$-reverse Carleson measures for Hardy spaces}\label{s_pq}
	
	Let $X$ be a quasi-Banach space of analytic functions in the unit disc and let $0< q< \infty$. 
	We say that $\mu\in M_+(\overline{\D})$ is a \emph{$q$-reverse Carleson measure} for $X$ if and only if there exists a constant $C > 0$ such that
	\begin{equation}
		\|f\|_{X}\leq C \|f\|_{L^q(\overline{\D},d\mu)} \label{q-RCM-def}
	\end{equation}
	for all $f\in X\cap C(\overline{\D})$.
	
	The main goal of the present section is to describe the $q$-reverse Carleson measures for the spaces $H^p(\mathbb D)$.
	
	\subsection{The balayage phenomenon}
	
	We start with the following simple observation: it turns out that if $X=H^p(\mathbb D)$, $0<p<\infty$,
	then only the part of the measure $\mu$ which lives on the boundary $\partial \D$ is `responsible' for the inequality \eqref{q-RCM-def}.
	
	\begin{lem}
		\label{general_balayage}
		Suppose that for a measure $\mu\in M_+(\overline{\D})$ there exists a constant $C > 0$ such that 
		\begin{equation}\label{q-RCM-Hardy}
			\|f\|_{H^p} \le C \|f\|_{L^q(\mu)}
		\end{equation}
		for every $f\in H^p(\D)\cap C(\overline{\D})$.
		Put $\nu=\mu|_{\partial\D}$. Then the inequality
		\begin{equation}
			\label{q-RCM-Hardy-balayaged}
			\|f\|_{H^p} \le C\|f\|_{L^q(\nu)}
		\end{equation}
		holds for every $f\in H^p(\D)\cap C(\overline{\D})$.
	\end{lem}
	
	\begin{proof}
		Let $f\in H^p(\D)\cap C(\overline{\D})$. Consider the following functions: $f_N(z)=z^N f(z)$. Apply the inequality \eqref{q-RCM-Hardy} to the functions $f_N$ to get
		\begin{equation}\label{RCM_fN}
			\|f_N\|_{H^p} \le C \|f_N\|_{L^q(\mu)}.
		\end{equation}
		Clearly, $\|f_N\|_{H^p} = \|f\|_{H^p}$, 
		since $|f_N(z)|=|f(z)|$ for $|z| = 1$. Besides that, for every $z\in\D$, we have $$|f_N(z)|=|z|^N |f(z)| \xrightarrow[N\to\infty]{}0.$$
		It is also clear that $|f_N(z)|\le \|f\|_{C(\overline{\D})}$ for all $z\in\overline{\D}$. Therefore, by the dominated convergence theorem, we have
		$$
		\lim_{N\to\infty} \int_\D |f_N(z)|^q\, d\mu(z) = 0.
		$$
		We now let $N$ tend to infinity in the inequality \eqref{RCM_fN} and obtain
		\begin{multline*}
			\|f\|_{H^p}\le C\lim_{N\to\infty} \Big( \int_\D |f_N(z)|^q\, d\mu(z) + \int_{\partial\D} |f_N(z)|^q\, d\mu(z) \Big)^{1/q}\\=C\Big(\int_{\partial\D} |f(z)|^q \, d\mu(z)\Big)^{1/q} = C\|f\|_{L^q(\nu)}, 
		\end{multline*}
		as required.
	\end{proof}
	
	\subsection{Main lemma}
	
	Now, we prove the main lemma, which helps us to describe the $q$-reverse Carleson measures for the spaces $H^p$. 
	This lemma says that we can substitute not only analytic functions into the inequality \eqref{q-RCM-Hardy-balayaged} (if we change there the $H^p$-norm to $L^p$-norm) but also arbitrary continuous functions.
	
	In what follows, it will be convenient for us to identify the unit circle $\T=\partial\D$ with the unit interval $[0,1)$.
	
	\begin{lem}\label{complex_analysis_goes_away}
		Suppose that $\mu$ is a $q$-reverse Carleson measure for $H^p$, where $p>0$ and $q>0$. Let $\nu=\mu|_{\partial\D}$, $\nu\in M_+(\T)$. Then there exists a constant $C > 0$ such that
		\begin{equation}\label{plugging_cont_f}
			\|f\|_{L^p(\T)}\leq C \|f\|_{L^q(\nu)}
		\end{equation}
		for every $f\in C(\T)$.
	\end{lem}
	
	\begin{proof}
		It is enough to prove the lemma for the trigonometric polynomials.
		Indeed, for any continuous function $f$, the inequality \eqref{plugging_cont_f} will then follow by the uniform approximation.
		
		Now, if $f$ is a trigonometric polynomial, then for a sufficiently large integer $N > 0$ the function $f_N(t)=e^{2\pi i N t} f(t)$ is an \emph{analytic} trigonometric polynomial: $\widehat{f}_N(n) = 0$ for all $n < 0$. Therefore, we can substitute this function into the inequality \eqref{q-RCM-Hardy-balayaged} and, clearly, the resulting estimate is equivalent to \eqref{plugging_cont_f}.
	\end{proof}
	
	\begin{rem}\label{reformulation}
		Obviously, if $\mu$ satisfies the hypothesis of Lemma~\ref{complex_analysis_goes_away}, 
		then it is a $q$-reverse Carleson measure for the space $H^p$. Hence, now we are left with the problem of characterization of measures $\nu\in M_+(\T)$ which satisfy the inequality \eqref{plugging_cont_f} for every $f \in C(\T)$.
	\end{rem}
	\subsection{The case $q < p$}
	
	\begin{theorem}
		If $0 < q < p$, then there are no $q$-reverse Carleson measures for the space $H^p$.
	\end{theorem}
	
	\begin{proof}
		Let $I\subset [0,1]$ be an interval of length $\eps > 0$ (we identify it with a subarc of $\partial\D$). 
		We apply inequality \eqref{plugging_cont_f} to the function $f=\mathbbm{1}_I$. 
		Observe that although this function is not continuous, it can clearly be pointwise approximated by a sequence of uniformly bounded continuous functions and we can apply the dominated convergence theorem.
		
		Inequality \eqref{plugging_cont_f} with $f=\mathbbm{1}_I$ yields
		$$
		\eps^{1/p} = m(I)^{1/p} \leq C \nu(I)^{1/q},
		$$
		and hence $\nu(I)\ge c \eps^{q/p}$. This inequality holds for every interval $I$ of length $\eps$. Now we consider approximately  $\eps^{-1}$ pairwise disjoint intervals of length $\eps$. Summing the corresponding inequalities for these intervals, we get $$\mu(\T) = \nu(\T)\ge c \eps^{\frac{q}{p}-1}.$$
		Since this inequality holds for arbitrarily small $\eps > 0$ and since $\frac{q}{p} - 1 < 0$, we can let $\eps$ tend to zero and conclude that the measure $\mu$ must be infinite. This contradiction finishes the proof.
	\end{proof}
	
	\subsection{The case $q > p$}
	\begin{theorem}
		Let $\mu\in M_+(\overline{\D})$ and let $\nu$ denote the restriction of $\mu$ to the boundary of the unit disc: $\nu=\mu|_{\partial\D}$. Let $\nu^a$ be the absolutely continuous part of $\nu$. Denote by $\beta$ the Radon--Nikodym derivative of $\nu^a$ with respect to Lebesgue measure on the unit circle. 
		If $0 < p < q$, then the measure $\mu$ is a $q$-reverse Carleson measure for the space $H^p$ if and only if 
		\begin{equation}\label{cond_beta}
			\frac{1}{\beta}\in L^{\frac{p}{q-p}}(\T).
		\end{equation}
	\end{theorem}
	\begin{proof}
		The `if' part is clear: if $\frac{1}{\beta}\in L^{\frac{p}{q-p}}(\T)$, then for any function $f\in H^p\cap C(\overline{\D})$,
		we apply H\"{o}lder's inequality with exponents $q/p$ and $(q/p)' = \frac{q}{q-p}$ and write
		\begin{multline*}
			\|f\|^p_{H^p} = \int_{\T} |f|^p \beta^{p/q} \beta^{-p/q}\, dm \le \Big(\int_\T |f|^q \beta\, dm\Big)^{p/q} \Big(\int_\T \beta^{-\frac{p}{q-p}}\, dm\Big)^{\frac{q-p}{q}} \\ = \|f\|_{L^q(\nu^a)}^p \Big\| \frac{1}{\beta} \Big\|_{L^{\frac{p}{q-p}}}^{p/q}\le \|f\|_{L^q(\mu)}^p \Big\| \frac{1}{\beta} \Big\|_{L^{\frac{p}{q-p}}}^{p/q}.
		\end{multline*}
		This means that $\mu$ is a $q$-reverse Carleson measure for the space $H^p$ and we are done.
		
		Suppose now that $\mu$ is a $q$-reverse Carleson measure for the space $H^p$.  
		Using Lemma~\ref{complex_analysis_goes_away}, we see that the property
		$$
		\|f\|_{L^p(\T)}\leq C \|f\|_{L^q(\nu)}
		$$
		holds for all continuous functions $f$. Consider an arbitrary nonnegative continuous function $f$ and denote $g = f^{p}$. We see that the following inequality holds for every nonnegative function $g\in C(\T)$:
		$$
		\Big( \int_{\T} g\, dm \Big)\le C \Big(\int_\T g^{q/p}\, d\nu \Big)^{p/q}. 
		$$
		Next, we apply the above inequality to the absolute value of an arbitrary function $g$, and conclude that the estimates
		$$
		\Big| \int_{\T} g\, dm \Big|\le \int_{\T} |g|\, dm\le C \Big(\int_\T |g|^{q/p}\, d\nu \Big)^{p/q}
		$$
		hold for all functions $g\in C(\T)$.
		Let $r = q/p > 1$. Then the functional
		$$
		\phi(g) =  \int_\T g\, dm 
		$$
		satisfies the bound $|\phi(g)|\le C \|g\|_{L^r(\nu)}$ for all functions $g\in C(\T)$ and hence it can be extended to a continuous linear functional on the whole space $L^r(\nu)$. This means that the measure $m$ is absolutely continuous with respect to $\nu$ and we have $m=\alpha\, d\nu$ with $\alpha\in L^{r'}(\nu)$. 
		The support of the function $\alpha$ is clearly contained in the support of the measure $\nu^a$. Hence, we can write $\nu^a = \frac{1}{\alpha}\, dm$. Let $\beta$ denote $1/\alpha$. 
		Then the condition $\int_\T \alpha^{r'}\, d\nu < \infty$ rewrites as follows:
		$$
		\int_{\T} \alpha^{r'}\, d\nu^a = \int_{\T}\alpha^{r'-1}\, dm = \int_{\T} \beta^{1-r'}\, dm < \infty.
		$$
		Since $r' = \frac{q}{q-p}$, we arrive at the required condition \eqref{cond_beta}.
	\end{proof}

	\section{Bloch, $\mathrm{BMOA}(\Dbb)$ and $H^\infty(\Dbb)$ spaces}\label{s_Bloch}
	Our main goal in this section is to  prove that there are no $q$-reverse Carleson measures for the Bloch space, for the space $\mathrm{BMOA}(\Dbb)$ and for the Hardy space $H^\infty(\Dbb).$ 
	
	\subsection{The Bloch and $\mathrm{BMOA}(\Dbb)$ spaces}
	Recall that the Bloch space $\mathcal{B}$ consists of those holomorphic in the open unit disc functions $f$ satisfying
	$$
	\sup_{z\in \Dbb}(1-|z|^2)|f'(z)|<\infty.
	$$
	The corresponding norm in the Bloch space is defined by
	$$
	\|f\|_{\mathcal{B}}=|f(0)|+\sup_{z\in \Dbb}(1-|z|^2)|f'(z)|.
	$$
	
	In turn, the space $\mathrm{BMOA}(\Dbb)$ consists of those  functions $f\in H^1(\mathbb D)$  that satisfy
	$$
	\|f\|_{\mathrm{BMOA}(\Dbb)}=|f(0)|+\sup_I \frac{1}{|I|}\int_{I}|f(e^{it})-f_I| dt < \infty,
	$$
	where the supremum is taken over all subarcs $I$ of the unit circle $\partial \Dbb$ and where we denote the average of $f$ 
	over $I$ as follows: $$f_I=\frac{1}{|I|}\int_I f.$$
	
	It is well known that $H^{\infty}(\Dbb)\subset \mathrm{BMOA}(\Dbb) \subset \mathcal{B}$.  Indeed, the former inclusion here is trivial, whereas the latter is an easy consequence of the Schwarz lemma. 
	Hence, once we manage to prove that, for a certain $q>0$, there exist no $q$-Carleson measures for the Bloch space, we will automatically prove the same for the spaces $\mathrm{BMOA}(\Dbb)$ and $H^\infty(\Dbb)$.
	
	Here is the main result of the current section.
	
	\begin{theorem}\label{Carl_Bloch}
		Let $q>0$. Then there are no $q$-reverse Carleson measures for the Bloch space.
	\end{theorem}
	
	We remark that in the second section of the current paper we have been considering $q$-reverse Carleson measures in certain $L^q$-based spaces. Since the Bloch $\mathcal B$ space is not tailored to a specific choice of $q$, it is hence natural to consider $q$-reverse Carleson measures in the space $\mathcal B$ for all $q>0$.
	
	\subsection{Proof of the main result}
	We will now concentrate on the proof of Theorem~\ref{Carl_Bloch}.

	Suppose that $\phi\in [-\pi,\pi]$ and that $n\in\mathbb N$ are both fixed. Define for $z\in \DDb$ an auxiliary function $f_n$ by the following formula:
	$$
	f_n(z)=\frac{ze^{in\phi}}{n}+\frac{z^2e^{i(n-1)\phi}}{n-1}+\ldots+\frac{z^{n-1}e^{i2\phi}}{2}+z^ne^{i\phi}
	.$$
	This function is nothing but a slightly modified $n$-th partial Fourier sum of the $n$-th analytic Fej\'er polynomial, see~\cite{pommerenke1974bloch}.
	
	First, observe that
	$$
	\|f_n\|_{\mathcal{B}}\gtrsim (1-|z_n|)|f'_n(z_n)|,
	$$
	where $z_n=e^{-1/n}e^{i\phi}$. Simple computations permit us to conclude that
	$$
	\|f_n\|_{\mathcal{B}}\gtrsim \frac{1}{n}\biggl(\frac{1}{n} +\frac{2e^{-1/n}}{n-1}+\ldots+\frac{(n-1)e^{-(n-2)/n}}{2}+ne^{-(n-1)/n}\biggr).
	$$
	Further, notice that the function $\psi(x)=x\cdot e^{-(x-1)/n}$ defined for real $x$ is increasing on the interval $[0,n]$. This allows us to  make use of the  Chebyshev  sum  inequality (see, for instance,~\cite{hardy1952inequalities}) and to infer the bound
	$$
	\|f_n\|_{\mathcal{B}}\gtrsim \frac{\log n}{n^2}\left(1+2e^{-1/n}+\ldots+(n-1)e^{-(n-2)/n}+ne^{-(n-1)/n}\right).
	$$
	From here we deduce that there exists a constant $C_0>0$ such that
	\begin{equation}
		\label{raz}
		\|f_n\|_{\mathcal{B}}\geq \frac{\log n}{C_0}
	\end{equation}
	for all $n \in\mathbb N$.
	
	Suppose now to the contrary that a finite measure $\mu$ defined in the closed unit disc is a $q$-reverse Carleson measure for a certain $q>0$. We will only give a detailed proof of  Theorem~\ref{Carl_Bloch} in the case where $q=1$, since the general case can be proven in a similar way. Indeed, the very same test function will work in the case of an arbitrary $q>0$. We hence infer that there exists a certain constant $C_1>0$ such that
	\begin{equation}
		\label{dva}
		\frac{\log n}{C_0}\leq C_1 \int_{\Dbb} |f_n(z)| d\mu(z) + C_1 \int_{\partial \Dbb} |f_n(z)| d\mu(z).
	\end{equation}
	
	We first concentrate ourselves on the estimate of the integral over the open unit disc above. We have
	\begin{equation}
		\begin{split}
			\label{tri}
			\int_{\Dbb} |f_n(z)| 
			&\, d\mu(z) \leq \int_{\Dbb} \left(\frac{|z|}{n}+\frac{|z|^2}{n-1}+\ldots+\frac{|z|^{n-1}}{2}+|z|^n\right) d\mu(z)\\
			&\leq \frac{1}{n}\left(\frac{1}{n}+\frac{1}{n-1}+\ldots+\frac{1}{2}+1\right)  \int_{\Dbb}  (|z|+|z|^2+\ldots+ |z|^n)d\mu(z)\\
			&\leq \frac{1+\log n}{n}  \int_{\Dbb}  \left(-1+\frac{1-|z|^{n+1}}{1-|z|}\right)d\mu(z),
		\end{split}
	\end{equation}
	where the second inequality above follows from the Chebyshev sum inequality. In the rest of the proof we suppose that $n$ is sufficiently large to guarantee that
	\begin{equation}
		\label{trispolovinoy}
		\frac{C_1 (1+\log n)}{n}   \int_{\Dbb} \left(-1+\frac{1-|z|^{n+1}}{1-|z|}\right)d\mu(z) \leq \frac{\log n}{4 C_0}.
	\end{equation}
	We indeed can choose such a number $n$ due to the dominated convergence theorem.
	
	We now proceed to the integral over the unit circle on the right-hand side of the inequality~\eqref{dva}. We infer the following chain of inequalities
	\begin{equation}
		\begin{split}
			\label{chetyre}
			\int_{\partial \Dbb} 
			&|f_n(z)| \,d\mu(z)=\int_{-\pi}^{\pi}\left|\frac{e^{in(\phi-t)}}{n}+\frac{e^{i(n-1)(\phi-t)}}{n-1}+\ldots+ e^{i(\phi-t)}\right| d\mu(e^{it}) \\
			&\leq \int_{\phi-a}^{\phi+a} (1+\log n) \; d\mu(e^{it}) + C_2' \int_{[-\pi,\pi)\backslash [\phi-a,\phi+a]} \log\left(\frac{10\pi}{|\phi-t|}\right) d\mu(e^{it})\\
			&\leq (1+\log n) \cdot \mu([\phi-a,\phi+a]) +C_2 \log\left(\frac{1}{a}\right) \mu(\partial \Dbb)
			\\
			&\le 2\log n \cdot \mu([\phi-a,\phi+a]) +C_2 \log\left(\frac{1}{a}\right) \mu(\partial \Dbb), 
		\end{split}
	\end{equation}
once $n$ is large enough and where $a\in (0,1/2)$ is to be chosen in a moment, $C_2, C_2'>0$ are absolute constants, and $[\phi-a,\phi+a]$ denotes the corresponding subarc of the unit circle. Here, the first inequality above is a direct consequence of the well-known bounds on the Dirichlet kernel and on its conjugate kernel, see, for instance,~\cite[Section 5.2]{zygmund2002trigonometric}.
	
	Let us now choose $a$. What we need is the following inequality:
	\begin{equation}
		\label{pyat}
		C_1\cdot C_2\cdot\mu(\partial \Dbb) \log\left(\frac{1}{a}\right)\leq \frac{\log n}{2 C_0}.
	\end{equation}
	This is certainly true if $a= n^{-1/(2C_0C_1C_2\mu(\partial \Dbb))}$, which is exactly the condition that we impose on $a$.
	
	Properties~\eqref{dva},~\eqref{tri},~\eqref{trispolovinoy},~\eqref{chetyre} and~\eqref{pyat} permit us to conclude that the estimate
	\begin{equation}
		\label{shest}
		\mu([\phi-a,\phi+a])\ge \frac{1}{8C_0C_1} - \frac{1}{8C_0C_1\log n} \ge \frac{1}{16C_0C_1}
	\end{equation}
	is verified by the measure $\mu$, provided that $n$ is sufficiently large. 
	
	Now, we consider $M\asymp a^{-1}$ points $\{\phi_k\}_{k=1}^M$ on the unit circle $\partial \Dbb$ that are vertices of a regular $M$-gon in such a way that the arcs $[\phi_k-a,\phi_k+a]$ are pairwise disjoint for $k$ from $1$ to $M$.  We apply the bound~\eqref{shest} to each of these intervals to deduce the following estimate $$\mu(\partial \Dbb)\gtrsim n^{1/(2C_0C_1C_2 \mu(\partial \Dbb))}.$$ Thus, we arrive at the contradiction with the fact that $\mu$ is finite, simply via letting $n$ tend to $\infty$. Hence, Theorem~\ref{Carl_Bloch} follows.

	\begin{rem}
		Note that the holomorphic Besov space $HB_{0}^{\infty,\infty}$ (see e.g.~\cite{OF99} or the sixth section of the current article for the definition of Besov spaces) coincides with the Bloch space. So, according to what has been proved in this section, there are no reverse Carleson measures there. On top of that, since for $p > 1$ one has $HB_{1/p}^{p,p} \subset  HB_0^{\infty,\infty}$, we see that there are no reverse Carleson measures in the spaces $HB_{1/p}^{p,p}$. 
		However, this fact has a more simple proof: see Section~\ref{s_besov} for the results concerning reverse Carleson measures for Besov spaces.
	\end{rem}

	\section{Triebel--Lizorkin spaces}\label{s_Triebel}
	\subsection{Main definitions and formulation of the result}
	Suppose that $s\ge 0$, $0 < p, q < \infty$. Then the holomorphic Triebel--Lizorkin space $HF_s^{p,q}$ in the unit disc is defined by the following seminorm:
	$$\|f\|^p_{F_{s}^{p,q}} = \int_0^1 \Big( \int_0^1 |f^{(m)}(re^{2\pi i t})|^q (1-r)^{(m-s)q - 1}\, dr \Big)^{p/q}\, dt,$$
	where $m$ is an arbitrary integer strictly greater than $s$ (the definitions for different choices of $m$'s are equivalent). As always, we can modify this definition for $q=\infty$:
	$$
	\|f\|^p_{F_{s}^{p,\infty}} = \int_0^1 \Big( \sup_{0\le r < 1} \big( |f^{(m)}(re^{2\pi i t})| (1-r)^{m-s} \big)\Big)^{p}\, dt.
	$$
	The parameter $s$ is sometimes called the smoothness of this space.
	
	The space $HF_s^{p,q}$ becomes a quasi-Banach space if we endow it with the quasinorm
	$$
	\|f\|_{HF_s^{p,q}} = \|f\|_{F_s^{p,q}} + \sum_{j=0}^{m-1} |f^{(j)}(0)|,
	$$
	see for instance~\cite{guliev1991b} or~\cite{OF99}. See also~\cite{tselishchev2020littlewood}
	for an equivalent norm within a certain range of parameters. 
	
	Our goal is to describe the reverse Carleson measures for Triebel--Lizorkin spaces, that is, to identify the measures $\mu\in M_+(\overline{\D})$ such that
	\begin{equation}\label{RCM_T-L}
		\|f\|_{HF_s^{p,q}} \leq C \|f\|_{L^p(\mu)}
	\end{equation}
	for every $f\in HF_s^{p,q}\cap C(\overline{\D})$.
	
	This problem is more general than its counterpart for Hardy spaces which was discussed in Section~\ref{s_hp}, since $HF_{0}^{p,2} = H^p$ (with equivalence of norms). This fact is well known and it follows from the Littlewood--Paley description of Hardy spaces. More generally, for any $s\ge 0$ we have $HF_{s}^{p,2} = H^p_s $, the latter space being the $p$-Hardy--Sobolev space of smoothness $s$ (see, for example,~\cite{manhas2018closures} for a definition).
	
	The main result of this section is the following theorem.
	
	\begin{theorem}\label{T-L_disc}
		Let $s\ge 0$, $0 < p < \infty$ and $0 < q \le\infty$. Then the following assertions hold.
		\begin{itemize}
			\item If $s > 0$, then  reverse Carleson measures for the space $HF_s^{p,q}$ do not exist.
			
			\item If $2\le q \le\infty$, then the measure $\mu$ is a reverse Carleson measure for the space $HF_0^{p,q}$ if and only if $\mu(I)\ge c|I|$ for every arc $I\subset \mathbb{\partial \D}$ and for some constant $c > 0$.
			
			\item If $0<q<2$, then the reverse Carleson measures for the space $HF_0^{p,q}$ do not exist.
		\end{itemize}
	\end{theorem}
	
	Observe that we consider only the case of nonnegative smoothness here. The main reason for this is that the case $s < 0$ is much better studied. In particular, for $s < 0$ the space $HF_s^{p,p}$ coincides with a certain weighted Bergman space, and the reverse Carleson measures for these spaces have been described, for example, in~\cite{Lu85}  and \cite{Lu2k}.  
	See also \cite{Lu88} for related results.
	
	The rest of the present section is devoted to the proof of  Theorem~\ref{T-L_disc}.
	
	\subsection{Positive smoothness}\label{sec:5.2}
	Let us start with the first (and the simplest) part of Theorem~\ref{T-L_disc}. 
	For $s_1 < s_2$, we have  $HF_{s_2}^{p,q}\subset HF_{s_1}^{p,q}$. Thus, 
	it is enough to consider the case $0 < s < 1$. Moreover, for every nonnegative $s$, we also have
	\begin{equation}\label{simple_embedd}
		HF_s^{p,q_1}\subset HF_s^{p,q_2} \quad \text{for} \quad q_1 < q_2,
	\end{equation} 
	and therefore we can consider only the case $q=\infty$.
	
	For the reader's convenience, we rewrite the definition of the seminorm in the Triebel--Lizorkin space for such parameters:
	\begin{equation}
		\|f\|_{F_s^{p,\infty}}^p = \int_0^1 \Big( \sup_{0\le r < 1} \big( |f'(re^{2\pi i t})| (1-r)^{1-s} \big)\Big)^{p}\, dt.
	\end{equation}
	Let us substitute $f_n(z)=z^n$ into this formula. If we prove that 
	\begin{equation}\label{to_prove_pos_smooth}
		\lim_{n\to\infty}\|f_n\|_{F_s^{p,\infty}} = \infty,
	\end{equation}
	then we are done, since on the other hand $\|f_n\|_{L^p(\mu)} \le \mu(\overline{\D})^{1/p}$, and this contradicts the inequality \eqref{RCM_T-L}.
	
	We have
	$$
	\|f_n\|_{F_s^{p,\infty}} = \sup_{0\le r < 1} (nr^{n-1}(1-r)^{1-s}).
	$$
	Now, we choose a particular value $r = \frac{n-1}{n-s} < 1$ and obtain
	\begin{multline}\label{norm_to_infty}
		\|f_n\|_{F_s^{p,\infty}} \gtrsim n\frac{(n-1)^{n-1}}{(n-s)^{n-s}} \ge \frac{(n-1)^n}{(n-s)^{n-s}} \\
		= \exp\bigl(n\log(n-1)-(n-s)\log(n-s)\bigr) \\ 
		= \exp\biggl(s\log(n-s)-n\log\Bigl(\frac{n-s}{n-1}\Bigr)\biggr).
	\end{multline}
	It is maybe worth noting that the value $r = \frac{n-1}{n-s}$ is obtained from the maximization of the quantity $nr^{n-1}(1-r)^{1-s}$.
	
	Thus, relation \eqref{to_prove_pos_smooth} follows, since the sequence $s\log(n-s)$ tends to infinity as $n\to \infty$, while the sequence $$n\log\biggl(\frac{n-s}{n-1}\biggr)\asymp \frac{n(1-s)}{n-1}$$ is bounded.
	
	\subsection{Zero smoothness, case $q \ge 2$}\label{sec:5.3} Now we pass to the proof of the second assertion in Theorem~\ref{T-L_disc}.
	Note that if a measure $\mu$ is such that $\mu(I)\gtrsim |I|$ for every arc $I\subset\partial\D$, then we have $$\|f\|_{HF_0^{p,q}} \lesssim \|f\|_{HF_0^{p,2}} \lesssim \|f\|_{H^p} \lesssim \|f\|_{L^p(\mu)}$$ for every function $f\in HF_0^{p,q}\cap C(\overline{\D})$.
	
	On the other hand, remark once again that in Section~\ref{s_hp} we proved the second assertion in Theorem~\ref{T-L_disc} for $q=2$ (since $HF_0^{p,2}=H^p$ with equivalence of norms). 
	Our characterization of reverse Carleson measures in this case is a consequence of the following fact. If 
	$$k_\lambda(z) = \frac{1}{1-\overline{\lambda}z}$$ and $lp > 1$, $l\in\mathbb{Z}$, 
	then $\|k_\lambda^l\|^p_{H^p}\asymp (1-|\lambda|)^{1-lp}$ for $|\lambda|$ close to $1$. 
	Hence, if we also manage to prove that 
	\begin{equation}\label{norm_kl_triebel}
		\|k_\lambda^l\|_{F_0^{p,q}}^p\asymp (1-|\lambda|)^{1-lp}
	\end{equation} 
	for all $|\lambda| > 1/2$ and $2\le q\le\infty$, then in a similar way we will readily derive the condition $\mu(I)\gtrsim |I|$ for every reverse Carleson measure $\mu$ for the space $HF_0^{p,q}$. 
	
	Now, observe that for any $q\ge 2$ we have $$\|k_\lambda^l\|_{HF_0^{p,q}}^p \lesssim \|k_\lambda^l\|_{H^p}^p \lesssim (1-|\lambda|)^{1-lp}.$$
	The last property is proved in Section 2 (see equation~\eqref{norm_klambda}). Thus, it remains only to show the inequality 
	\begin{equation}\label{key_norm_est}
		\|k_\lambda^l\|_{F_0^{p,q}}^p \gtrsim (1-|\lambda|)^{1-lp}.
	\end{equation}
	Owing to the property \eqref{simple_embedd}, we hence infer that it is enough to prove this inequality for $q=\infty$.
	
	Recall that
	$$
	\|f\|_{F_0^{p,\infty}}^p = \int_0^1 \Big( \sup_{0\le r < 1} \big( |f'(re^{2\pi i t})| (1-r) \big)\Big)^{p}\, dt.
	$$
	Let $f=k_\lambda^l$. For $|\lambda| > 1/2$, we obtain
	\begin{equation}\label{element_est_1}
		|f'(r e^{2\pi i t})|=l |k_\lambda'(re^{2\pi i t})|\cdot |k_\lambda (re^{2\pi i t})|^{l-1} \gtrsim \frac{1}{|1-\overline{\lambda}re^{2\pi i t}|^{l+1}}.
	\end{equation}
	Without loss of generality, we can assume that $\lambda \in \mathbb{R}$, $\lambda > 1/2$. 
	Put $\alpha = 1-\lambda$. A direct computation shows that for $0 < t < \alpha/4$ the following inequality holds:
	\begin{equation}\label{element_est_2}
		|1-(1-\alpha)^2 e^{2\pi i t}| \leq 6 \alpha.
	\end{equation}
	Indeed, we have
	\begin{multline*}
		|1-(1-\alpha)^2 e^{2\pi i t}| \leq |1-(1-\alpha)^2 \cos 2\pi t| + |(1-\alpha)^2 \sin 2\pi t| \\ \leq |1-\cos 2\pi t| + |(2\alpha - \alpha^2)|\cdot |\cos 2\pi t| + |\sin 2\pi t| \leq 2\pi^2 t^2 + 2\alpha + |2\pi t| \leq 6\alpha.
	\end{multline*}
	For $f=k_\lambda^l$, substituting $r=1-\alpha$ instead of the supremum and then using  \eqref{element_est_1} and \eqref{element_est_2}, we obtain
	\begin{multline*}
		\|f\|_{F_0^{p,\infty}}^p \ge \int_0^{\alpha/4} \Big( \alpha |f'((1-\alpha)e^{2\pi i t})|  \Big)^{p}\, dt \\ \gtrsim \int_0^{\alpha/4} \Big( \frac{\alpha}{|1-(1-\alpha)^2 e^{2\pi i t}|^{l+1}} \Big)^p\, dt   \gtrsim \alpha^{1-lp},
	\end{multline*}
	and therefore estimate \eqref{key_norm_est} is proved and we are done.
	
	\subsection{Zero smoothness, case $q < 2$} Now we prove the final assertion of Theorem~\ref{T-L_disc}. 
	Observe that due to \eqref{simple_embedd}, it is enough to consider only the case $q\ge 1$. 
	In fact, we prove a more general statement: the following lemma is an extension of a result of Rudin \cite{Ru55}. It seems to be interesting in its own right.
	
	\begin{lem}\label{q-var}
		For any $1\le q < 2$, there exists a function $f\in H^\infty$ such that 
		\begin{equation}\label{q-var-inf}
			\int_0^1 |f'(r e^{2\pi i t})|^q (1-r)^{q-1} dr = \infty \ \ \text{for a.e.}\  t\in (0,1). 
		\end{equation}
	\end{lem}
	
	Clearly, if $f$ is a function provided by Lemma~\ref{q-var} then $f\not\in HF_{0}^{p,q}$ for any $p > 0$. Therefore, the inequality \eqref{RCM_T-L} cannot hold (formally, we can consider the functions $f_\rho(z)=f(\rho z)$, plug them into the inequality \eqref{RCM_T-L} and take the limit as $\rho\to 1-$). 
	
	For $q=1$, Lemma~\ref{q-var} was proved by Rudin \cite{Ru55}: the proof is based on the estimates for lacunary Fourier series proved by Zygmund \cite{Zyg44}. We follow a similar route but instead of using the results by Zygmund we use their generalization from \cite{Wat50}. 
	
	The proof of Lemma~\ref{q-var} consists of two steps.
	
	\subsubsection{Step 1: construction of a function in $H^2$} First, we construct a function $g\in H^2$ such that 
	\begin{equation}\label{inf_q_var}
		\int_0^1 |g'(r e^{2\pi i t})|^q (1-r)^{q-1} dr = \infty \ \ \text{for every}\  t\in (0,1). 
	\end{equation}
	Consider the following function:
	$$
	g(z) = \sum_{n=1}^\infty \frac{z^{2^n}}{n^{1/q}}.
	$$
	Since $q < 2$, we have
	$$
	\|g\|^2_{H^2} = \sum_{n=1}^\infty \frac{1}{n^{2/q}} < \infty,
	$$
	and hence $g\in H^2$.
	
	Now we are in position to use \cite[Theorem 2]{Wat50}. Let us state this theorem here.
	
	\begin{theorem*} 
		(D. Waterman.) Suppose that $0 < \lambda_1 < \lambda_2 < \ldots$ and $\lambda_{k+1}/\lambda_k > c > 1$ for every $k\ge 1$. Consider the Dirichlet series $h(s) = \sum_{n=1}^\infty a_n e^{-\lambda_n s}$. Then for every $q > 1$ the following inequality holds:
		$$
		\sum_{n=1}^\infty |a_n|^q \le C \int_0^\infty (1-e^{-s})^{q-1} |h'(s)|^q\, ds.
		$$
	\end{theorem*}
	
	Since $\sum_{n=1}^\infty \frac{1}{n} = \infty$, making the substitution $r=e^{-s}$,  we see that the above theorem indeed implies that relation \eqref{inf_q_var} holds.
	
	\subsubsection{Step 2: construction of a bounded function} Denote the expression on the left-hand side of equation \eqref{q-var-inf} by $V_q(f; t)$. Now, we construct a function $f\in H^\infty$ such that $V_q(f; t) = \infty$ for a.e. $t$. It can be done in exactly the same way as in \cite{Ru55}; however, we prefer to briefly describe this construction here for completeness.
	
	The function $g$ which was obtained in the previous step can be written as $g=h_1/h_2$, where $h_1, h_2\in H^\infty(\D)$. We put $A_j = \{t\in (0,1): V_q(h_j; t) = \infty\}$, $j= 1, 2$, and $A=A_1\cup A_2$. Since $V_q(g;t) = \infty$ for every $t$, $|A| = 1$. Indeed, if $h_2(e^{2\pi i t})\neq 0$, which is true for a.e.\ $t$, then condition \eqref{inf_q_var} implies that at least one of the expressions $V_q(h_1; t)$ and $V_q(h_2;t)$ must be infinite.
	
	Now, for every $w\in\mathbb{C}$, denote $f_w = h_1 + wh_2$ and $E_w = \{t\in A : V_q(f_w; t) < \infty\}$. If we prove that the sets $E_w$ are pairwise disjoint, then we are done, since it will imply that $|E_w| = 0$ for all but countably many values of $w$.
	
	Suppose that $t\in A_1\cap E_w$. Then for every $v\neq w$ we can write $wf_v = vf_w + (w-v)h_1$ and hence $t\not \in E_v$. Analogously, if $t\in A_2\cap E_w$, then for $v\neq w$, using the formula $f_v = f_w + (v-w)h_2$, we conclude that $t\not\in E_v$. Therefore, the above argument shows that we can put $f = h_1 + wh_2$ for some $w\in\mathbb{C}$, and we are done.

	\section{Besov spaces}\label{s_besov}
	
	In this section, we consider yet another important scale of function spaces, that is, holomorphic Besov spaces. Using the recent results from~\cite{baranov2024analytic} and one of the methods from the previous section (which goes back to~\cite{HMNO14}), 
	that is, substitution of the functions $k_\lambda$, we describe the reverse Carleson measures for Besov spaces. We consider only the case $p\ge 1$ here. It is an interesting question whether Theorem~\ref{t_Besov} below can be extended to all values $p > 0$.
	
	\subsection{Main definitions and formulation of the result}
	
	Recall the definition of Besov spaces of analytic functions (see e.g.~\cite{baranov2024analytic}). Let $0 < p, q < \infty$ and $s\ge 0$. The holomorphic Besov space $HB_{s}^{p,q} = HB_s^{p,q}(\D)$ consists of functions $f \in \hol(\mathbb D)$ such that
	$$\int_0^1\left(\|f^{(m)}(re^{2\pi i(\cdot)})\|_{L^p(\mathbb T)} (1-r)^{m-s-1/q}\right)^q dr< \infty,$$
	where $m$  is any integer strictly greater than $s$ (the definitions for different choices of such $m$'s are equivalent).
	A quasinorm in this space is given by the following expression
	$$
	\|f\|_{HB^{p,q}_s}:= \sum_{j=0}^{m-1} |f^{(j)}(0)| +\left(\int_0^1\left(\|f^{(m)}(re^{2\pi i(\cdot)})\|_{L^p(\mathbb T)} (1-r)^{m-s-1/q}\right)^q dr\right)^{\frac{1}{q}}.
	$$
	This definition can be extended to infinite values of $p$ or $q$ in the standard way. It is worth noting that $HB_0^{\infty, \infty} = \mathcal{B}$, the Bloch space that we have considered in Section~\ref{s_Bloch}.
	The problem of describing the reverse Carleson measures for Besov spaces is mentioned, e.g. at the end of the paper~\cite{hartmann2021dominating}.
	
	Our goal is to describe all reverse Carleson measures for Besov spaces, i.e., all measures $\mu\in M_+(\overline{\D})$ satisfying
	\begin{equation}
		\label{besov1}
		\|f\|_{ HB^{p,q}_s} \leq C \|f\|_{L^p(\mu)},
	\end{equation}
	for all $f\in HB^{p,q}_s\cap C(\overline{\mathbb D})$.
	As always, we refer to the parameter $s$ above as the smoothness of the corresponding Besov space.
	
	\begin{theorem}\label{t_Besov}
		Let $s\ge 0$, $1 \le p < \infty$ and $0 < q \le\infty$. Then the following assertions hold.
		\begin{enumerate}
			\item If $s > 0$, then there exist no reverse Carleson measures for the space $HB_s^{p,q}$.
			
			\item If $ p \ge 2$ and $q \ge p$, then the measure $\mu$ is a reverse Carleson measure for the space $HB_0^{p,q}$ if and only if $\mu(I)\ge c|I|$ for every arc $I\subset \mathbb{\partial \D}$ and for some constant $c > 0$.
			
			\item If $p\ge 2$ and $q < p$ then there exist no reverse Carleson measures for the space $HB_0^{p,q}$.
			
			\item If $p < 2$ and $q\ge 2$ then the measure $\mu$ is a reverse Carleson measure for the space $HB_0^{p,q}$ if and only if $\mu(I)\ge c|I|$ for every arc $I\subset \mathbb{\partial \D}$ and for some constant $c > 0$.
			
			\item If $p < 2$ and $ q < 2$ then there exist no reverse Carleson measures for the space $HB_0^{p,q}$.
		\end{enumerate}
	\end{theorem}
	
	Now we briefly present the proof of the above theorem. We constantly use the  embedding
	\begin{equation}\label{simple_embedd_bes}
		HB_s^{p,q_1}\subset HB_s^{p,q_2} \quad \text{for} \quad q_1 < q_2.
	\end{equation}
	Besides that, we use the fact that Besov and Triebel--Lizorkin spaces are related in the following way. 
	
	\begin{pr}\label{emb_tr_bes}
		Besov and Triebel--Lizorkin spaces satisfy the following embeddings:
		\begin{itemize}
			\item If $p\le q$ then $HF_0^{p,q}\subset HB_0^{p,q}$;
			
			\item If $p\ge q$ then $HB_0^{p,q}\subset HF_0^{p,q}$
		\end{itemize}
	\end{pr}
	
	For a proof of this proposition see e.g.~\cite[Theorem 4.1]{OF99}, where the corresponding embeddings are presented in the more general context of spaces of functions in the unit ball of $\mathbb{C}^d$.

	\subsection{Proof of Theorem~\ref{t_Besov}}
	
	\subsubsection{Positive smoothness}
	
	The case $s > 0$ can be treated similarly to that of Triebel--Lizorkin spaces: it is enough to apply the inequality~\eqref{besov1} to the functions $f_n(z) = z^n$. 
By \eqref{simple_embedd_bes}, it suffices to consider the case $q=\infty$.
A simple computation shows that for $0 < s < 1$
	$$
	\|f_n\|_{HB_s^{p,\infty}} = \sup_{0\le r < 1} (nr^{n-1}(1-r)^{1-s}) =\|f_n\|_{HF_s^{p,\infty}},
	$$
	and that this last quantity above tends to $\infty$ (along the lines of  Section~\ref{s_Triebel}). Since the functions $f_n$ are uniformly bounded, this implies the first claim of Theorem~\ref{t_Besov}.
	
	\subsubsection{Case $q\ge p\ge 2$}
	
	Recall the definition of the functions 
	$$
	k_\lambda(z) = \frac{1}{1-\overline{\lambda}z}.
	$$
	We will consider these functions for $1/2 < |\lambda| < 1$. Due to Proposition~\ref{emb_tr_bes} and the results of Section~\ref{s_Triebel} (see equation~\eqref{norm_kl_triebel}), for $q\ge p\ge 2$, we have
	$$
	\|k_\lambda\|^p_{HB_0^{p,q}}\lesssim \|k_\lambda\|^p_{HF_0^{p,q}}\asymp (1-|\lambda|)^{1-p}.
	$$
	Therefore, in order to prove the second assertion of Theorem~\ref{t_Besov}, we only need to prove the estimate
	\begin{equation}\label{key_est_bes}
		\|k_\lambda\|^p_{HB_0^{p,q}}\gtrsim (1-|\lambda|)^{1-p},
	\end{equation}
	since the rest of the proof follows in a standard way (as in Section~\ref{s_hp} above). Due to the relation~\eqref{simple_embedd_bes}, it is enough to prove the above inequality~\eqref{key_est_bes}  for $q = \infty$.
	
	Using Lemma~\ref{very_simple}, for $|\lambda| > 1/2$, we write
	\begin{multline*}
		\|k_\lambda\|_{HB_0^{p,\infty}} = \sup_{0 < r< 1} \Big\{ \Big( \int_0^1 |k_\lambda'(re^{2\pi i t})|^p\, dt \Big)^{1/p} (1-r) \Big\}\\ \asymp \sup_{0 < r < 1} \Big((1-|\lambda|r)^{\frac{1-2p}{p}}(1-r)\Big) = \sup_{0 < r < 1}\Big( \frac{1-r}{(1-|\lambda|r)^{2-\frac{1}{p}}}\Big).
	\end{multline*}
	The substitution  $r = |\lambda|$ in the last supremum above gives us the required estimate~\eqref{key_est_bes}, and we are done.
	
	\subsubsection{Case $p\ge 2$, $q < p$}
	This case follows from the results of~\cite{baranov2024analytic}. For $z\in \mathbb C$, we denote
	$$
	B_n(z):= \prod_{k=1}^n \frac{z^{2^k}-a}{1-az^{2^k}},
	$$
	where $a:=1-1/n$. It can be easily seen that the functions $B_n$ are Blaschke products, simply via factorization.  
	
	We note further that for $2 \leq p$ and $1\leq q \leq p <\infty$, it holds that
	\begin{equation}
		\label{blashkesozvezdoj}
		\|B_n\|_{HB^{p,q}_0}\gtrsim n^{\frac{1}{q}-\frac{1}{p}},
	\end{equation}
	see~\cite[proof of Proposition 10]{baranov2024analytic}. Since the Blaschke products $B_n$ are uniformly bounded, we can use them in the inequality~\eqref{besov1}. Letting $n$ tend to $\infty$ we  get a contradiction which implies the third statement of Theorem~\ref{t_Besov} for $q\ge 1$. The general case $q > 0$ follows from~\eqref{simple_embedd_bes}.
	
	\subsubsection{Case $p < 2\le q$}
	If $p > 1$ then the proof of the fourth statement of our theorem is completely similar to the proof of the second one above; we omit it. For $p=1$ we use the test function $k_\lambda^2$ in the place of $k_\lambda$, as in Section~\ref{sec:5.3}.
	
	\subsubsection{Case $p, q < 2$}
	
	It has also been noticed in~\cite[Proposition 10]{baranov2024analytic} (see   the end of Section 2 in~\cite{barankayum} for a proof) that for $(p,q)\in [1,2]^2$ there exists a Blaschke product $\widetilde{B}_n$ such that
	$$
	\|\widetilde{B}_n\|_{HB^{p,q}_0}\gtrsim (\log n)^{\frac{1}{q}-\frac{1}{2}}.
	$$
	
	Similarly to the third case, it also means that there are no reverse Carleson measures for the spaces $HB_0^{p,q}$, once $1\leq q<2$, $1\leq p\leq 2$, and the general case $0 < q < 2$ once again follows from the relation~\eqref{simple_embedd_bes}.

\section{Higher dimensions}\label{s_highdim}

The problem of description of reverse Carleson measures for spaces of holomorphic functions in the unit ball of $\mathbb{C}^d$, $d\ge 2$, becomes much more difficult. The analogue of Theorem~\ref{t_Carl_Hp_onedim} still holds in this situation, however, we are not able to prove other results in the same generality as in the rest of the paper. We illustrate it by considering Triebel--Lizorkin spaces, see Theorem~\ref{t_Triebel_multidim} below. Some other related questions are listed as open problems in the next section.

\subsection{Main definitions and reverse Carleson measures for Hardy spaces}

Let $\bd$ denote the open unit ball of $\ccd$, $d\ge 1$
and let $\si=\sid$ denote the Lebesgue measure on the boundary unit sphere $\spd$. We denote by $\hol(\bd)$ the space of holomorphic functions in the ball $\bd$.

For $0<p<\infty$ the Hardy space $H^p=H^p(\bd)$
consists of functions $f\in \hol(\bd)$ satisfying
\[
\|f\|_{H^p}^p = \sup_{0<r<1} \int_{\spd} |f(r\zeta)|^p\, d\si(\zeta) < \infty.
\]
As usual we identify the Hardy space $H^p(\bd)$, $p>0$ with the space
$H^p(\spd)$ which consists of the corresponding boundary values.
We refer the reader to the monographs \cite{Ru80} or \cite{Z05} for more information on the Hardy spaces $\hp(\bd)$.

Consider $\MM(\bd)$, the space of finite, positive Borel measures on the ball $\bd$.

For an arbitrary dimension $d\ge 1$ the following sets play a role similar to that of  the subarcs of the circle:
\[
Q = Q(\za, \delta) =\{\xi\in \spd: |1-\langle \za, \xi\rangle | \le \delta\},
\quad \za\in\spd,\ \delta>0.
\]
Remark that $Q$ is a closed ball with respect to the non-isotropic metric
\[
\rho(\za, \xi) = |1-\langle \za, \xi\rangle |^{\frac{1}{2}}, \quad \za, \xi\in\spd.
\]
For $r>0$ and $Q = Q(\za, \delta)$ we put $rQ = Q(\za, r\delta)$ by definition.

The theorem about Carleson measures in the unit ball (see for instance \cite[Chapter~5]{Z05})
yields that the bound
\begin{equation}\label{e_Carl_n}
\| f \|_{L^p(B_d, \mu)} \le C \|f\|_\hp,\quad f \in\hp,
\end{equation}
is equivalent to the following condition on the measure $\mu$:
\begin{equation}\label{e_Carl_geom_n}
\mu(S_Q) \le C\si(Q)\quad\textrm{for all non-isotropic balls}\ Q\subset\spd,
\end{equation}
where
\[
S_Q = \left\{z \in\bd:\ 1-\delta \le |z| \le 1,\ \frac{z}{|z|} \in  Q \right\}
\]
is a standard Carleson window for $Q= Q(\cdot, \delta)$.

We are now interested in the following reverse Carleson inequality:
\begin{equation}\label{e_revrese}
\|f \|_\hp \le C \|f\|_{L^p(\bbd,\mu)}, \quad f\in  C(\bbd) \cap H^p(\bd),\ 0 < p < \infty,
\end{equation}
for $\mu\in\MM(\bbd)$.
 
Fix a number $l\in\Nbb$ such that $pl >1$. For $w \in \bd$ consider the functions
\[
k_w(z) = \frac{1}{(1-\langle z, w\rangle)^{d}},\quad z\in\bd,
\]
and the corresponding auxiliary normalized versions defined by
\[
K_w = \frac{k_w^l}{\|k_w^l\|_{\hp}}.
\]

By \cite[Proposition~1.4.10]{Ru80},
\[
\frac{(1-|w|^2)^{\frac{d-p\ell d}{p}}}{C}\leq \|k_w^l\|_p \leq C(1-|w|^2)^{\frac{d-p\ell d}{p}}
\]
for a certain constant $C>1$ independent of $w$.

The generalization of Theorem~\ref{t_Carl_Hp_onedim} now reads as follows.

\begin{theorem}\label{t_Carl_Hp}
Let $0 < p < \infty$ and let $\mu \in \MM(\bbd)$, $d\ge 1$.
The following assertions are equivalent.
\begin{enumerate}
	\item[(i)] There exists a constant $C > 0$ such that
	\[
	\int_{\bbd}
	|f|^p\,d\mu \ge C \|f\|_\hp^p\quad\textrm{for all}\ f \in C(\bbd) \cap H^p.
	\]
	\item[(ii)] There exists a constant $C > 0$ such that
	\[
	\int_{\bbd}
	|K_w|^p\,d\mu \ge C \quad\textrm{for all}\ w\in\bd.
	\]
	\item[(iii)] There exists a constant $C > 0$ such that
	\[
	\mu(Q) \ge C \si(\qball)\quad\textrm{for all non-isotropic balls}\ Q \subset \partial\bd.
	\]
\end{enumerate}
\end{theorem}

For $p > 1$ and $l = 1$ this theorem is proved in~\cite{Dou24}. The general case can be deduced similarly to the one-dimensional case (see discussion after Theorem~\ref{t_Carl_Hp_onedim} above). On the other hand, one can deduce this theorem from the multidimensional analogues of Lemmas~\ref{very_simple} and~\ref{simple}. Since these computations are purely technical generalizations of the computations in Section~\ref{s_hp} of our paper (and they are similar to the computations in~\cite{Dou24}), we omit them.

\begin{rem}
The properties from Theorem~\ref{t_Carl_Hp} are equivalent to the following:
there exists a constant $C> 0$ such that
\begin{equation}\label{e_iv}
\mu(S_\qball) \ge C \si(\qball)\quad\textrm{for all non-isotropic balls}\ Q \subset \partial\bd.
\end{equation}
The proof of this fact can be found in~\cite[Proposition 1]{Dou24}.
\end{rem}

\subsection{Triebel--Lizorkin spaces in the unit ball}

Now we consider Triebel--Lizor\-kin spaces of holomorphic functions in the unit ball of $\mathbb{C}^d$. Recall that for a function $f\in\hol(B_d)$ its radial derivative  $\rad f$ is defined by the following formula:
\[
\rad f(z) = \sum_{k=1}^{d} z_k \frac{\partial f}{\partial z_k}(z), \quad z \in \bd.
\] 
Thus, if $f(z) =\sum_\alpha c_\alpha z^\alpha\in \hol(\bd)$, then for $s\in \Rbb$ we define
\[
(I+\rad)^s f(z) =\sum_\alpha (1+|\alpha|)^s c_\alpha z^\alpha.
\]

By definition, the holomorphic Triebel--Lizorkin space $HF_{s}^{p,q} = HF_{s}^{p,q}(B_d)$  consists of functions $f\in \hol(\bd)$ satisfying
\[
\|f\|_{HF^{pq}_s}^p = \int_{\spd} \left( \int_{0}^1 |(I + \rad)^{[s]^+} f(r\zeta)|^q 
(1-r^2)^{([s]^+-s)q-1}\, dr \right)^{\frac{p}{q}}\, d\sigma(\zeta)
<\infty,
\]
where $[s]^+$ is the integral part of the number $s+1$. It is not difficult to see that for $d = 1$ this definition is equivalent to the definition of Triebel--Lizorkin spaces from Section~\ref{s_Triebel}.

As usual, this definition can be generalized to the case $q = \infty$. Notice that, similarly to the one-dimensional case, a norm in the space $HF^{p,2}_0$ is defined as the $L^p$-norm of the Littlewood--Paley $g$-function,
so the space $HF^{p,2}_0$ coincides with the Hardy space $H^p$. It is also well known that for $s > 0$ the space $HF^{p,2}_s$ coincides with the holomorphic Hardy--Sobolev space $H_s^p$ in the unit ball of $\mathbb{C}^d$. 
We refer the reader to the paper~\cite{OF99} for more information on holomorphic Triebel--Lizorkin spaces. In particular, the following statement will be useful for us.

\begin{pr}[see {\cite{BB89}, \cite[Theorem~4.1]{OF99}}]\label{p_embed}
Let $s\in \Rbb$ and let $0<p<\infty$. Then
\[
H F^{p,q_1}_s \subset H F^{p,q_2}_s, \quad 0< q_1 \le q_2 \le \infty.
\]
\end{pr}

Our goal in the present section is to prove the $d$-dimensional counterpart of Theorem~\ref{T-L_disc}. 
We say that $\mu\in M_+(\overline{B}_d)$ is a reverse Carleson measure for the space $H F^{pq}_s$ if there exists a constant $C > 0$ such that
\begin{equation}\label{rev_Carl_Triebel_multidim}
\|f\|_{HF^{pq}_s}\le C\|f\|_{L^p(\mu)}
\end{equation}
for every $f\in HF^{pq}_s\cap C(\overline{B}_d)$.

\begin{theorem}\label{t_Triebel_multidim}
Let $s\ge 0$, $0 < p < \infty$ and $0 < q \le\infty$. Then the following assertions hold.
\begin{itemize}
\item If $s > 0$, then there exist no reverse Carleson measures for the space $HF_s^{p,q}$.

\item If $2\le q \le\infty$, then the measure $\mu$ is a reverse Carleson measure for the space $HF_0^{p,q}$ if and only if $\mu(Q)\ge c \sigma(Q)$ for every non-isotropic ball $Q\subset \mathbb{\partial B_d}$ and for some constant $c > 0$.

\item If $0 < q \le 1$, then there exist no reverse Carleson measures for the space $HF_0^{p,q}$.
\end{itemize}    
\end{theorem}

Clearly, the range of exponents $1 < q < 2$ is missing here. The main reason for this is that we do not know if the analogue of Lemma~\ref{q-var} holds in several dimensions. See the next section, where we discuss certain related open problems for further details. 

\subsubsection{Positive smoothness}
The proof of this part for $d = 1$ in Section~\ref{s_Triebel} was very simple: we used the functions $f_n(z)=z^n$. In our present situation the role of these functions will be played by the following version of Ryll-Wojtaszczyk polynomials \cite{RW83}. Arguing as in Section~\ref{sec:5.2} we may assume that $0 < s < 1$.

\begin{pr}[see e.g. \cite{A86}]
Let $d\ge 2$ and $0<p<\infty$.
There exist $J \in \Nbb$ and $\delta>0$, and for each $n\in\Nbb$
there exist holomorphic homogeneous polynomials $W_{n j}$, $1\le j \le J$, on $\ccd$
of degree $n$, such that
\begin{align}
\sum_{j=1}^J |W_{n j}(\za)|^p &\le 1, \quad \za\in \spd, \label{e_RWup}\\
\max_{1\le j \le J} |W_{n j}(\za)| &\ge \delta>0,   \quad \za\in \spd.  \label{e_RWlow}
\end{align}
\end{pr}

Assume that $\mu$ is a reverse Carleson measure for $HF_{s}^{p,\infty}$ with $s > 0$.
On the one hand, \eqref{e_RWup} implies that
\begin{equation}\label{e_Pnj_mu}
\sum_{j=1}^{J} \|W_{n j}\|^p_{L^p(\mu)}\le \mu (\overline{B}_d).
\end{equation}

On the other hand, \eqref{e_RWlow} guarantees that
\[
\sum_{j=1}^{J} \|W_{n j}\|_{HF_s^{p, \infty}} \ge \delta\sup_{0\le r < 1} ((n+1) r^{n} (1-r)^{1-s}),
\]
since $W_{n j}$ is a holomorphic homogeneous polynomial of degree $n$.
By the estimates obtained for $d=1$ (see the end of Subsection~5.2 and relation~\eqref{norm_to_infty}), for $s > 0$ we have
\[
\sup_{0\le r < 1} ((n+1) r^{n} (1-r)^{1-s}) \to \infty \ \textrm{as}\ n\to \infty.
\]
Therefore,
\[
\sum_{j=1}^{J} \|W_{n j}\|_{HF_s^{p, \infty}} \to \infty \ \textrm{as}\ n\to \infty.
\]
The above property combined with \eqref{e_Pnj_mu} implies that the reverse Carleson measures 
do not exist for the $HF_s^{p, \infty}$, $s>0$, and hence, due to Proposition~\ref{p_embed}, also for all spaces  $HF_s^{p, q}$ with $s > 0$ and $0 < p, q < \infty$.

\subsubsection{Zero smoothness, case $q\ge 2$}

Similarly to the one-dimensional case, it suffices to prove that
\begin{equation}\label{e_kwlow}
\|k_w^\ell\|^p_{HF_{0}^{p, q}} \gtrsim (1-|w|)^{d-p\ell d}
\end{equation}
for $q=\infty$.
Recall that
\[
\|f\|^p_{HF_{0}^{p, \infty}} = \int_{\spd} \left(\sup_{0\le r < 1} (|(\rad +I) f|(1-r))\right)^p\, d\sid(\za).
\]
Therefore,
\[
\|k_w^\ell\|^p_{HF_{0}^{p, \infty}} \gtrsim 
\int_{\spd} \left(\sup_{0\le r < 1} \frac{1-r}{|1-\langle \za, r w\rangle|^{\ell d +1}} \right)^p\, d\sid(\za).
\]
Putting $r=|w|$, we obtain
\[
\|k_w^\ell\|^p_{HF_{0}^{p, \infty}} \gtrsim
\int_{\spd} \left( \frac{1-|w|}{|1-\langle \za, |w| w\rangle|^{\ell d +1}} \right)^p\, d\sid(\za)
\gtrsim (1-|w|)^p (1-|w|)^{-p\ell d - p + d}
\]
by \cite[Proposition~1.4.10]{Ru80}.
Hence, the inequality \eqref{e_kwlow} holds.

\subsubsection{Zero smoothness, case $q \le 1$}
Due to Proposition~\ref{p_embed}, it is enough to prove that there do not exist reverse Carleson measures for the space $HF_0^{p,1}$. This fact can be deduced from the following result by Rudin, which is in turn a generalization of the result from the paper~\cite{Ru55} to the high-dimensional setting.
\begin{pr}[see {\cite[Theorem~17.9]{Ru86}}]
There exists a function $F\in H^\infty(B_d)\cap C(\overline{B}_d)$ such that
$$
\int_0^1 |(I+\rad) F(r\zeta)|\, dr = \infty
$$
for a.e. $\zeta\in\spd$.
\end{pr}
Clearly, such function $F$ has  infinite norm in the Triebel--Lizorkin space $HF_0^{p,1}$. Hence, the inequality~\eqref{rev_Carl_Triebel_multidim} cannot hold (formally, we can consider the functions $F_\rho(z)=F(\rho z)$, plug them into the inequality \eqref{rev_Carl_Triebel_multidim} and take the limit as $\rho\to 1-$).

\section{Concluding remarks and open problems}\label{s_remarks}

\subsection{Besov and Triebel--Lizorkin spaces}

We have not tried to cover all possible exponents for Besov spaces like we did when we considered Triebel--Lizorkin spaces. Indeed, in Theorem~\ref{t_Besov} we have not considered the case $p < 1$. However, if $p < 1$ and $q \ge 2$ one can prove that a measure $\mu$ is a reverse Carleson measure for the space $HB_0^{p,q}$ if and only if $\mu(I)\ge C|I|$ for every arc $I\subset \mathbb{\partial \D}$ in a standard way: by the application of the functions $k_\lambda^l$ with an integer $l$ such that $pl > 1$.

On the other hand, if $0 < q \le p < 1$, then there are no reverse Carleson measures for the space $HB_0^{p,q}$: it can be deduced from the embedding $HB_0^{p,q}\subset HF_0^{p,q}$ (see~\cite[Theorem~4.1]{OF99}) and the fact that in this case there are no reverse Carleson measures for the space $HF_0^{p,q}$ due to Theorem~\ref{T-L_disc}.

Hence only the following case remains open.

\begin{problem}
Is it true that if $0 < p < 1$ and $p < q < 2$ then there are no reverse Carleson measures for the space $HB_0^{p,q}$?
\end{problem}
For $p\ge 1$ in the proof of the corresponding statement we used the powerful results from~\cite{baranov2024analytic} about the existence of certain Blaschke products.
However, for our purposes it would be enough to construct a function with similar properties that belongs to the space $H^\infty(\mathbb D)$. 
It would therefore be interesting to obtain a more elementary proof, which could  work for all exponents $p$ and solve the above problem.

We also note that in Sections~\ref{s_Triebel} and~\ref{s_besov} we considered only the 
reverse Carleson inequality in the form $\|f\|_{X_0^{p,q}}\lesssim \|f\|_{L^p(\mu)}$, where $X$ stands for either $HF$ or $HB$, Triebel--Lizorkin or Besov spaces,
respectively. The main motivation was that the inequalities of such kind are similar to those from Theorem~\ref{t_Carl_Hp_onedim}. However, it seems natural for us to consider the following more general question.

\begin{problem}
Is it possible to describe the measures $\mu\in M_+(\overline{\mathbb{D}})$ such that the inequality
$$
\|f\|_{X_0^{p,q}}\lesssim \|f\|_{L^s(\mu)}
$$
holds for all functions $f\in X_0^{p,q}\cap C(\overline{\D})$? Here, once again, $X$ stands for either $HF$ or $HB$.
\end{problem}
This question does seem technically complicated. On the other hand, it is possible that a positive answer to this question entails a technique which is similar to our new method applied in Section~\ref{s_pq}. It would hence be interesting to obtain such general results.

\subsection{Several dimensions}

As we have mentioned above, the things become much more difficult in higher dimensions, even if we consider reverse Carleson measures for function spaces in the unit ball of $\mathbb{C}^d$. The reason for this is that already the problems of \emph{existence} of analytic functions in the unit ball with certain specific properties are highly nontrivial (probably the most well-known example is the problem of existence of inner functions in the unit ball). Therefore, it would be interesting to obtain the generalizations of all our results from Sections 3--6 to the higher dimensional setting. Let us mention one particular problem, whose positive solution would lead to the proof of Theorem~\ref{t_Triebel_multidim} for the missing range of exponents $1 < q < 2$.

\begin{problem}
Let $1 < q < 2$. Does there exist a function $F\in H^\infty(B_d)$ such that the property
$$
\int_0^1 |(I+\rad) F(r\zeta)|^q (1-r)^{q-1}\, dr = \infty
$$
holds for a.e. $\zeta\in\spd$?
\end{problem}
This would be an interesting generalization of our Lemma~\ref{q-var} to higher dimensions. For $q = 1$ this problem, as we have mentioned, was solved by Rudin~\cite{Ru86}, but, unlike the one-dimensional case, in several dimensions the generalization to the case of all exponents $q < 2$ seems to be complicated.

\bibliographystyle{plain}
\bibliography{ReCarl}
\end{document}